\documentclass[a4paper,12pt]{article}
\usepackage{amssymb,amsmath}
\usepackage[dvips]{color}
\usepackage[dvips]{graphics}
\newtheorem{theorem}{Theorem}[section]
\newtheorem{lemma}{Lemma}[section]
\newtheorem{remark}{Remark}[section]
\newtheorem{example}{Example}

\newcommand{\eps}{\varepsilon_M}

\usepackage{algorithmic}
\usepackage{algorithm}

\begin{document}
\title{On  Constructing  Orthogonal Generalized Doubly Stochastic Matrices}
\author{Gianluca Oderda\thanks{Ersel Asset Management SGR S.p.A., Piazza Solferino, 11, 10121 Torino, Italy, e-mail:  
Gianluca.Oderda@ersel.it}, 
Alicja Smoktunowicz\thanks{Faculty of Mathematics and Information Science, Warsaw University of Technology, 
Koszykowa 75, 00-662 Warsaw, Poland, e-mail: A.Smoktunowicz@mini.pw.edu.pl}
and Ryszard Kozera\thanks{Faculty of Applied Informatics and Mathematics, Warsaw University of Life Sciences - SGGW,
Nowoursynowska str. 159, 02-776 Warsaw, Poland and 
Department of Computer Science and Software Engineering, The University of Western Australia, 35 Stirling Highway, 
Crawley, WA 6009, Perth, Australia, e-mail:  ryszard.kozera@gmail.com}}

\maketitle

\medskip

\begin{abstract}
A real quadratic matrix is  generalized doubly stochastic (g.d.s.) if all of its row sums and column sums equal one. 
We propose numerically stable methods for generating such matrices having possibly  
orthogonality property or/and satisfying Yang-Baxter equation (YBE). 
Additionally, an inverse eigenvalue problem for finding orthogonal 
generalized doubly stochastic matrices with prescribed eigenvalues is solved here.
The tests performed in \textsl{MATLAB} illustrate our proposed algorithms and 
demonstrate their useful numerical properties.
\end{abstract}

 \noindent {\bf AMS Subj. Classification:} {15B10, 15B51, 65F25, 65F15.}  


\section{Introduction}
\label{Introduction}
We propose efficient algorithms  for constructing generalized doubly stochastic matrix  $A  \in  \mathbb R^{n \times n}$. 
Recall that $A$  is a generalized doubly stochastic matrix (g.d.s.) if all of its row sums and column sums equal one.
Let $I_n$ denote the $n \times n$  identity matrix and $e=(1,1, \ldots, 1)^T=\sum_{i=1}^n e_i \in  \mathbb R^{n}$, where
$\{e_i\}_{i=1}^n$ forms a canonical basis in $\mathbb R^{n}$.
The set 
\[
\mathcal{A}_n= \{A  \in  \mathbb R^{n \times n}:  Ae=e, \quad  A^Te=e\}
\]
of all such g.d.s. matrices is investigated in this paper. Noticeably, the class of g.d.s. matrices
$\mathcal{A}_n$  includes a thinner subset of all doubly stochastic matrices (bistochastic) -  see \cite{Horn}, pp. 526-529.  
However, in contrast to  the latter,  a  generalized doubly stochastic matrix does not necessarily permit only non-negative entries. 

Let $\mathcal{B}_n$ define the space of orthogonal generalized doubly  stochastic matrices determined by the following condition:
\[
\mathcal{B}_n= \{Q  \in  \mathcal{A}_n:  Q^TQ=I_n\}.
\]
 
Some applications of doubly stochastic matrices or g.d.s. matrices are outlined in \cite{brualdi}-\cite{Glunt2}. 
More specifically, in economy,  the orthogonal  generalized doubly stochastic matrices 
permit to map a space of original quantities (asset prices) into a space of transformed asset prices. 

Recall that if  $A  \in  \mathbb R^{n \times n}$ is bistochastic and orthogonal then $A$ is actually a permutation matrix (see e.g. \cite{Horn}). 
The situation is different for orthogonal generalized doubly stochastic matrices. Indeed, as simple inspection reveals
\[
A=\frac{1}{3} \left(
\begin{array}{ccc}
-1 &  2 & 2 \\
2 &  2 & -1 \\
2 &  -1 & 2
\end{array}
\right)
\]
forms an orthogonal generalized doubly stochastic matrix evidently not yielding a permutation matrix. 

\medskip

We address now the question of how to construct generalized doubly stochastic matrices and orthogonal g.d.s. matrices.
Let us denote by $\mathcal{Q}_n$  the set of all orthogonal matrices of size $n$:
\[
\mathcal{Q}_n=\{Q  \in  \mathbb R^{n \times n}: Q^T Q=I_n\},
\] 
and define
\[
\mathcal{U}_n=\{Q\in  \mathcal{Q}_n:  q_1=Qe_1=\frac{1}{\sqrt{ n}} e\}.
\]

The  following theorem will  be suitable later for the construction of some 
herein proposed algorithms. 
\begin{theorem}\label{thm2}
Given any $Q \in \mathcal{U}_n$ and  any $X \in  \mathbb R^{(n-1) \times (n-1)}$. Define   
\begin{equation}\label{B}
B=\left(
\begin{array}{cc}
1 &  0^T \\
0 &  X
\end{array}
\right).
\end{equation}
Then $A=QBQ^T$ is a generalized doubly stochastic matrix. 

On the other hand, 
if $A \in  \mathbb R^{n \times n}$ is a g.d.s.  matrix and  $Q \in \mathcal{U}_n$ then for $B=Q^TAQ$ we have (\ref{B})
 for some $X \in  \mathbb R^{(n-1) \times (n-1)}$. 

Moreover, $A$ is orthogonal if and only if $X$ defined in (\ref{B}) is orthogonal.
\end{theorem}
\begin{proof}  First, observe that from (\ref{B}) it follows that  $Be_1=e_1$ and $B^T e_1=e_1$.
Since  $Q \in \mathcal{U}_n$  we have $Q^T e=\sqrt{n} e_1$,  and so
\[
A e= QB(Q^Te)=\sqrt{n} Q(Be_1)=\sqrt{n} Qe_1= e.
\]
Similarly, 
\[
A^T e= QB^T(Q^Te)=\sqrt{n} Q(B^Te_1)=\sqrt{n} Qe_1= e.
\]

The proof for $B=Q^TAQ$ may be handled analogously.

Clearly, $A=QBQ^T$ is orthogonal for any orthogonal matrix $B$.
\end{proof}

\begin{remark}\label{glunt}
Theorem \ref{thm2} is  a slight reformulation of  the result established in \cite{Glunt2},  for $q_n=Q e_n$ instead of $q_1=Qe_1$.
\end{remark}

Note  that if $A_1,  A_2 \in \mathcal{A}_n $  then  $(A_1+A_2)/2 \in \mathcal{A}_n $ and $A_1 A_2 \in \mathcal{A}_n$.
Clearly, if $A_1,  A_2 \in \mathcal{B}_n $  then also  $A_1 A_2 \in \mathcal{B}_n$. 
Visibly, the latter renders various possible schemes for the derivation of the orthogonal doubly stochastic matrices. 

This paper focuses on constructing  orthogonal  generalized doubly  stochastic matrices with additional 
special properties enforced. More specifically,
in Section \ref{algorithms} some new algorithms for generating matrix $Q \in \mathcal{U}_n$ 
using the Householder QR decomposition (see e.g. \cite{Golub})
are proposed. We also describe a method for constructing   $A \in \mathcal{B}_n$ and propose the new 
algorithms for computing orthogonal generalized doubly stochastic matrices with prescribed eigenvalues. 
At the end of Section \ref{algorithms} a new scheme for constructing
orthogonal  generalized doubly stochastic matrices satisfying the Yang-Baxter equation (YBE) is also given.
Section \ref{numexp} includes numerical examples all implemented in \textsl{MATLAB} illustrating the new methods
introduced in this work.  Finally, the Appendix annotating this paper
includes the respective codes in \textsl{ MATLAB} for all algorithms in question.

\medskip

\section{Algorithms}
\label{algorithms}
The Algorithms $1-6$  for constructing orthogonal generalized doubly stochastic matrices are proposed and 
discussed below.
The respective \textsl{ MATLAB} codes of implemented algorithms  are attached in the Appendix.

\subsection{Construction of a symmetric $A \in \mathcal{B}_3$}
The aim is now to find a symmetric matrix $A \in \mathcal{B}_3$ in the following form:
\begin{equation}\label{matrixA}
A=\left(
\begin{array}{ccc}
x &  y & z \\
y &  z & x \\
z &  x & y
\end{array}
\right), 
\end{equation}
which also satisfies
\begin{equation}\label{eqs1}
x+y+z=1,
\end{equation}
and meets the orthogonality conditions:
\begin{equation}\label{eqs2}
x^2+y^2+z^2=1,
\end{equation}
\begin{equation}\label{eqs3}
x y +yz + xz=0.
\end{equation}

Clearly,  the equation (\ref{eqs3}) follows from  (\ref{eqs1})-(\ref{eqs2}) due to:
\[
1=(x+y+z)^2= x^2+y^2+  z^2+ 2(xy + yz+xz).
\]
Furthermore by (\ref{eqs1}) we obtain:
\begin{equation}\label{eqs4}
x+y=1-z.
\end{equation}
Hence $x^2+y^2 +z^2=(x+y)^2-2 xy +z^2=(1-z)^2-2xy+z^2$, which together with (\ref{eqs2}) yields:
\begin{equation}\label{eqs5}
xy=z(z-1). 
\end{equation}

For a given real number $z$ the solution $x$ of (\ref{eqs4})-(\ref{eqs5}) should satisfy the quadratic equation
$x^2-x(1-z)-z(1-z)=0$.  Since $\Delta= (1-z)(1+3z)$ we conclude that $x$ remains 
real if and only if  $z  \in  [-1/3,1]=I$. 
In this case we have two real solutions: $x_1=(1-z+\sqrt{\Delta})/2$ and $x_2=(1-z- \sqrt{\Delta})/2$.
Consequently, for $i=1,2$ two pairs of real solutions 
$(x_i, y_i)$ satisfying   (\ref{eqs4}) and  (\ref{eqs5}) can be now found according to the
procedures specified below (Algorithm 1 for $x=x_1$ and Algorithm 1a for $x=x_2$). 
However, the choice $x=x_2$ in Algorithm 1a leads to severe loss 
of accuracy of the computed  result once $z$ gets very close to $0$. 
Indeed, here two nearly equal numerator's numbers  $1-z$ and $\sqrt{\Delta}$ are then subtracted yielding
an undesirable effect of {\em ``nearly zero cancellation''}. 
In contrast, for $z\in I$ the Algorithm 1 does not bear such
computational deficiency adding merely two positive numbers in its numerator, respectively.
For more details see \cite{Nick}, Sec. 1.8. Solving a Quadratic equation, pp. 10-12. 
The comparison between two methods demonstrating the above mentioned cancellation pitfall 
is given later in Example \ref{example1a}.

\bigskip

{\bf Algorithm 1.}  {\em Construction of  $A \in \mathcal{B}_3$ of the  form (\ref{matrixA}).}

\bigskip

Choose first an arbitrary $z  \in  [-1/3,1]$. 

The algorithm consists of the following steps:
\begin{itemize}
\item If $z=1$ then $A=\left(
\begin{array}{ccc}
0 &  0 & 1 \\
0 &  1 & 0 \\
1 &  0 & 0
\end{array}
\right)$. 
\item If $-1/3 \leq z < 1$ then compute
\begin{itemize}
\item $\Delta= (1-z)(1+3z)$,
\item $x_1=(1-z+\sqrt{\Delta})/2$,
\item $y_1=-z(1-z)/x_1$,
\item $A=\left(
\begin{array}{ccc}
x_1 &  y_1 & z \\
y_1 &  z & x_1 \\
z &  x_1 & y_1
\end{array}
\right)$. 
\end{itemize}
\end{itemize}

\bigskip

{\bf Algorithm 1a (unstable for $z \approx 0$).} {\em Construction of  $A \in \mathcal{B}_3$. }

\bigskip

Choose first an arbitrary $z  \in  [-1/3,1]$. 

The algorithm consists of the following steps:
\begin{itemize}
\item If $z=1$ then $A=\left(
\begin{array}{ccc}
0 &  0 & 1 \\
0 &  1 & 0 \\
1 &  0 & 0
\end{array}
\right)$. 
\item If $z=0$ then 
$A=\left(
\begin{array}{ccc}
0 &  1 & 0 \\
1 &  0 & 0 \\
0 &  0 & 1
\end{array}
\right)$ .
\item  if  $z \neq 0$ and $-1/3 \leq z < 1$  then compute
\begin{itemize}
\item $\Delta= (1-z)(1+3z)$,
\item $x_2=(1-z-\sqrt{\Delta})/2$,
\item $y_2=-z(1-z)/x_2$,
\item $A=\left(
\begin{array}{ccc}
x_2 &  y_2 & z \\
y_2 &  z & x_2 \\
z &  x_2 & y_2
\end{array}
\right)$. 
\end{itemize}
\end{itemize}

\subsection{Construction of  $Q \in \mathcal{U}_n$ by Householder QR method}\label{Un}

In this subsection, we resort to the Householder method for computing the QR factorization of a given matrix 
$X  \in  \mathbb R^{n \times n}$. Recall that in
\textsl{ MATLAB},  the statement $[Q,R] = qr(X)$  decomposes $X$ into
an upper triangular matrix $R \in  \mathbb R^{n \times n}$ and orthogonal  matrix 
$Q  \in  \mathbb R^{n \times n}$ so that $X= Q R$. This method uses a suitably chosen sequence of Householder transformations. 
The reason for selecting the Householder method 
instead of the others including e.g.
Gram-Schmidt orthogonalization methods, is that the Householder QR decomposition
is {\em unconditionally  stable} (see, e.g. \cite{Nick}, Chapter $18$).

Recall that a Householder transformation (Householder reflector) is a matrix of the form
\[
H=I_n-\frac{2}{z^Tz} z z^T, \quad 0 \neq z \in \mathbb R^{n}.
\]

Note that $H$ is symmetric and orthogonal.  Householder matrices are very useful while
introducing zeros into vectors to transform matrices into simpler forms (e.g. triangular, bidiagonal etc.).

For example, if $z=e+\sqrt{n} e_1$ is taken then 
\[
He=(I_n-\frac{2}{z^Tz} z z^T) e=e- \frac{2}{z^Tz} z (z^T e)=-\sqrt{n} e_1.
\]

Similarly $He_1=- \frac{1}{\sqrt{ n}} e$ and therefore $Q=-H \in \mathcal{U}_n$.

\bigskip

In this paper a different  algorithm (Algorithm 2) for computing $Q \in \mathcal{U}_n$
based on Householder QR decomposition is proposed.
It enables to generate a vast class of orthogonal matrices with the first column equal to $\frac{1}{\sqrt{ n}} e$.

\medskip

{\bf Algorithm 2.} {\em Construction of  $Q \in \mathcal{U}_n$.}

\medskip

Let  $X=(x_1,x_2, \ldots, x_n)  \in  \mathbb R^{n \times n}$ be an arbitrary quadratic  matrix, with
each $x_i\in  \mathbb R^{n}$.

The subsequent steps read as:
\begin{itemize}
\item $q=\frac{1}{\sqrt{ n}} e$,
\item $\hat X= (q,x_2, \ldots, x_n)$,
\item $\hat X= \hat Q \hat R$ (Householder QR factorization),
\item $Q=-\hat Q$. 
\end{itemize}

\medskip

\begin{remark}\label{remark1}  Note that if $Q  \in \mathcal{U}_n $ and $Z \in  \mathcal{Q}_n$ 
is an arbitrary orthogonal matrix such that $Ze_1=e_1$ then   $Q Z \in \mathcal{U}_n$, and therefore there
 are many other  choices to  create the matrix  $Q \in \mathcal{U}_n$.
\end{remark}

\bigskip

\subsection{General method for constructing $A \in \mathcal{B}_n$}

Note also that Theorem \ref{thm2} permits to establish a general method for generating orthogonal g.d.s. matrix. 
Indeed the following scheme accomplishes such task:

\medskip

{\bf Algorithm 3.} {\em Construction of  $A \in \mathcal{B}_n$.}

\medskip

Take first arbitrary  $Q \in \mathcal{U}_n$ and  $W  \in \mathcal{Q}_{n-1}$.

The algorithm is determined now by two steps:
\begin{itemize}
\item 
$B=\left(
\begin{array}{cc}
1 &  0^T \\
0 &  W
\end{array}
\right),
$
\item $A=QBQ^T$.
\end{itemize}

\medskip

At this point, we remark that at the preliminary step $Q$  can be generated
by Algorithm $2$ and $W$ can be determined upon applying  Householder QR decomposition.

\medskip 

\subsection{Orthogonal generalized doubly stochastic matrices with prescribed  eigenvalues}

This subsection focuses on constructing the orthogonal generalized doubly stochastic matrix  
with prescribed eigenvalues. 
In doing so, a real Schur decomposition of orthogonal matrices is applied.
More specifically, recall a well-known result (Theorem 7.4.1  in \cite{Golub}):
\begin{theorem}{\bf (Real Schur Decomposition)}\label{thm3}
If $A \in  \mathbb R^{n \times n}$, then there exists an orthogonal $W \in  \mathbb R^{n \times n}$  and 
$R \in  \mathbb R^{n \times n}$ such that $A=WRW^T$, where 
\begin{equation}\label{rschur}
R= 
\left(
\begin{array}{cccc}
R_{11} & R_{12}  & \dots  & R_{1s}\\
0 & R_{22}  & \dots & R_{2s} \\
\vdots & \vdots  & \ddots &  \vdots \\
0  & 0  & \dots  &   R_{ss}
\end{array}
\right),
\end{equation}
and each $R_{kk}$ is either a $1$-by-$1$ matrix or a $2$-by-$2$ matrix having complex conjugate eigenvalues. 
\end{theorem}

\begin{lemma}\label{lemat1} If $A$ in Theorem \ref{thm3} is additionally
orthogonal then $R$ in (\ref{rschur}) is also orthogonal and hence $R$ is a block diagonal matrix
\begin{equation}\label{ortschur}
R = 
\left(
\begin{array}{cccc}
R_{11} & 0   & \dots  & 0\\
0 & R_{22}  & \dots & 0 \\
\vdots & \vdots & \ddots &  \vdots \\
0  & 0  & \dots   &  R_{ss}
\end{array}
\right),
\end{equation}
with each $R_{kk}$ forming either $\pm 1$  or a $2$-by-$2$ real matrix having complex conjugate eigenvalues 
$z_k=c_k+i s_k$ and $\bar z_k= c_k - i s_k$, where $c_k^2+s_k^2=1$.   
\end{lemma}

\medskip

We apply now Lemma \ref{lemat1} to generate a special form (\ref{ortschur}) of orthogonal g.d.s. matrices.  

\noindent

{\bf Algorithm 4.} {\em  Construction of  $A \in \mathcal{B}_n$ with prescribed eigenvalues.}

\medskip

{\bf Input:} 
\begin{itemize}
\item $r$- the number of the eigenvalues of $A$  equal to $1$, $r \ge 1$,
\item $p$- the number of the  eigenvalues of $A$  equal to $-1$,  $p \ge 1$,
\item given  $z=(z_1, z_2, \ldots, z_m)^T  \in  \mathbb C^{m}$- the vector of the eigenvalues of $A$, $m \ge 1$,
\item given arbitrary  $Q \in \mathcal{U}_n$, where $n=r+p+2m$. 
\end{itemize}

{\bf Output:} $A \in \mathcal{B}_n$ having  the eigenvalues $\pm 1$, and  $z_k$, $\bar z_k$ for $k=1, \ldots, m$.

The subsequent steps of the algorithm obey the following  pattern:
\begin{itemize}
\item Find $c_k$ and $s_k$ such that  $z_k=c_k+i s_k$ ($c_k$ is the real part and $s_k$ is the imaginary part of $z_k$), 
for $k=1, \ldots, m$,  
\item compute the rotation matrices $M_k$, for $k=1, 2, \ldots, m$
\begin{equation}\label{Rk}
R_k = 
\left(
\begin{array}{cc}
 c_k & s_k\\
-s_k & c_k
\end{array}
\right),
\end{equation}
\item create  a block diagonal matrix $R (2m \times 2m)$ 
\begin{equation}\label{M1}
R = 
\left(
\begin{array}{cccc}
R_{1} & 0   & \dots  & 0\\
0 & R_{2}  & \dots & 0 \\
\vdots & \vdots & \ddots &  \vdots \\
0  & 0  & \dots   &  R_{m}
\end{array}
\right),
\end{equation}
\item form the matrix $B$ according to:
\[
B = 
\left(
\begin{array}{ccc}
 I_r & 0 & 0\\
 0 & -I_p & 0\\
 0 & 0 & R\\
\end{array}
\right),
\]
\item compute $A=QBQ^T$.
\end{itemize}

\medskip

\begin{remark}\label{remark3}
Note that  $R_k$ defined by
(\ref{Rk}) is an orthogonal matrix  with the eigenvalues equal to $c_k+i s_k$ and $c_k-i s_k$.  
Clearly, one can extend Algorithm 4 to the special cases of $p=0$ or $m=0$. It is omitted here for the
sake of brevity. Noticeably, the case of $r=0$ in Algorithm 4 is excluded.
\end{remark}

\medskip

\subsection{Construction of  $A \in \mathcal{B}_n$ satisfying the Yang-Baxter equation}

Recall that matrix $A \in \mathbb R^\mathrm{n^2 \times n^2}$ satisfies the  Yang-Baxter Equation (YBE) if 
\begin{equation}\label{Eqs1}
(A \otimes I_n)    (I_n \otimes A)   (A \otimes I_n)   = (I_n  \otimes A)    (A \otimes I_n)   (I_n \otimes A), 
\end{equation}
where $X \otimes Y$ is the Kronecker product (tensor product) of the matrices $X$ and $Y$:
$X \otimes  Y= (x_{i,j} Y)$. That is, the  Kronecker product $X \otimes  Y$ is a block matrix whose $(i,j)$ blocks are $x_{i,j} Y$.

The Yang-Baxter equation has been extensively studied due to its application in many fields of mathematics
or quantum information science - for detailed applications see  e.g. \cite{Agatka}.
Solutions of the Yang-Baxter equation have many interesting properties.
Of particular importance to this work is the following  theorem  (see  \cite{Agatka}): 
\begin{theorem}\label{thm4}
If  $A \in \mathbb R^{n^2 \times n^2}$  satisfies the Yang-Baxter equation  (\ref{Eqs1}) and $P \in \mathbb R^{n\times n}$  is arbitrary non-singular matrix, then 
$\hat X=  (P \otimes P)A(P \otimes P)^{-1}$  also  satisfies the  Yang-Baxter equation (\ref{Eqs1}).
\end{theorem}

Based on the latter the efficient  algorithm (see \cite{Agatka}),  for generating special
 solutions of  the  Yang-Baxter equation (\ref{Eqs1}) can be now formulated.
 
\medskip

\noindent

{\bf Algorithm 5.} {\em Construction of  $A=A(d) \in \mathbb R^{n^2 \times n^2}$ satisfying the  Yang-Baxter equation.}

Select an arbitrary  $d=(d_1, d_2, \ldots, d_{n^2})^T \in \mathbb R^{n^2}$.

The algorithm obeys the following pattern:
\begin{itemize}
\item  Form $n^2$-by-$n^2$ matrix $S$:
\begin{equation}\label{Sn}
S=\left(\begin{array}{cccccc}
1 & 2 & 3 & \dots & n-1 & n\\
n+1 &  n+2 & n+3 & \dots & 2n-1& 2n \\
2n+1 & 2n+2 & 2n+3 & \dots & 3n-1 & 3n\\
\vdots & \vdots & \vdots & \vdots & \vdots &  \vdots \\
\vdots & \vdots & \vdots & \vdots & \vdots &  \vdots \\
\vdots & \vdots & \vdots & \vdots & \vdots &  \vdots \\
(n-1)n+1 &  (n-1)n+2 & (n-1)n+3 & \dots & n^2-1&   n^2
\end{array}\right),
\end{equation}
\item take $p=(p_1, p_2, \ldots, p_{n^2}) = (s_1^T, s_2^T, \ldots, s_{n^2}^T)$, where $s_j$ denotes the  $j$th column of $S$,
\item  set  $A=(d_{p_1} \, e_{p_1}, d_{p_2}\,  e_{p_2},  \ldots,  d_{p_{n^2}} \, e_{p_{n^2}})$,
\item then define $A=DP$, where $D=diag(d_1,  d_2, \ldots, d_{n^2})$ and $P=(e_{p_1}, e_{p_2},  \ldots, e_{p_{n^2}})$ is a permutation matrix.
\end{itemize}

\medskip

\begin{remark}
Note that the matrix  $X$ generated by Algorithm $5$  satisfies
 $Xe_1=d_1 e_1$ and $X^T e_1=d_1 e_1$ for arbitrary $d=(d_1, d_2, \ldots, d_{n^2})^T \in \mathbb R^{n^2}$.
In particular, upon taking  $n=2$ and $d=(d_1,d_2,d_3,d_4)^T$ we arrive at:
\[
X=\left(
\begin{array}{cccc}
d_1 &  0 & 0  & 0\\
0 &  0 & d_2  & 0\\
0  &  d_3 & 0  & 0\\
0 &  0 & 0  & d_4
\end{array}
\right).
\] 

More detailed information can be found in \cite{Agatka}.
\end{remark}

In order to generate the orthogonal solutions to the YBE we prove now the following:

\begin{theorem}\label{thm5}
Let $P \in \mathcal{U}_n$.  Assume that $B \in \mathbb R^{n^2 \times n^2}$ is an  orthogonal matrix  satisfying  the Yang-Baxter equation: 
\[
(B \otimes I_n)    (I_n \otimes B)   (B \otimes I_n)   = (I_n  \otimes B)  (B \otimes I_n)   (I_n \otimes B),
\]
with the additional conditions  $Be_1= e_1$ and $B^T e_1=e_1$.  

Define $Q= P \otimes P$ and $A=Q B Q^{T}$.  Then $Q \in   \mathcal{U}_{n^2}$  and $A \in  \mathcal{B}_{n^2}$ is orthogonal and  satisfies the  Yang-Baxter equation (\ref{Eqs1}).
\end{theorem}

\begin{proof}
Observe that $Q$ is orthogonal since 
\[
Q^T Q=  (P^T \otimes P^T)  (P \otimes P)=  (P^TP) \otimes (P^TP)=I_n \otimes I_n=I_{n^2}.
\]

We shall verify now that $Q e_1=\frac{1}{n} \bar e$, where $\bar e= (1,1, \ldots, 1)^T \in  \mathbb R^{n^2}$. 
Clearly,  $\bar e = e \otimes e$, where $e= (1,1, \ldots, 1)^T \in  \mathbb R^{n}$.  

Since  $P \in \mathcal{U}_n$ we have $P \hat e_1= \frac{1}{\sqrt{ n}} e$, where 
$\hat e_1=(1,0, \ldots, 0)^T \in  \mathbb R^{n}$. Hence we obtain $P^T e= \sqrt{ n} \hat e_1$. 

Exploiting now the standard properties of the Kronecker product yields:
\[
Q^T \bar e=  (P^T \otimes P^T) (e \otimes e)= (P^T e) \otimes  (P^T e)= n (\hat e_1) \otimes (\hat e_1)= n e_1,
\]
and so finally $Q e_1= \frac{1}{\sqrt{n^2}} \bar e$. The proof is complete. 
\end{proof}

\medskip 
Having established Theorem \ref{thm5}, we pass now to the formulation of the last algorithm.
\medskip

\noindent

{\bf Algorithm 6.} {\em Construction of  orthogonal generalized doubly stochastic matrix $A \in \mathbb R^{n^2 \times n^2}$ satisfying the  Yang-Baxter equation (\ref{Eqs1}).}

\medskip

Let $P \in \mathcal{U}_n$ and $B \in \mathbb R^{n^2 \times n^2}$ form an arbitrary matrix 
satisfying the assumptions of Theorem \ref{thm5}.

The algorithm splits into two steps:
\begin{itemize}
\item  $Q= P \otimes P$,
\item $A=QBQ^T$.
\end{itemize}

\medskip

In order to initialize the above procedure, the matrix $P$ is obtainable from Algorithm $2$, whereas
$B$ is computable with the aid of Algorithm $5$, where $d_1=1$ and 
$d_2, \ldots, d_{n^2}$ are arbitrary parameters satisfying  $|d_i|=1$, for all $i$. 

\medskip 

\section{Numerical Experiments}
\label{numexp}

The final section of this paper reports on  the results of the numerical experiments examining
the computational properties of Algorithms 1-6. 
All tests are performed in \textsl{MATLAB} version 8 \textsl{.4.0.150421 (R2014b)}, with machine precision $\eps \approx 2.2 \cdot 10^{-16}$.

We report on the following statistics for a given matrix $A$:
\begin{itemize}
\item ${err}_{orth}= \|I_n- A^T A \|_2$ (the orthogonality error),\label{orth}
\item ${err}_{rows}= \|A e-e\|_2$ (the error in the row sums),
\item ${err}_{columns}= \|A^T e-e\|_2$ (the error in the column sums).
\end{itemize}

Here $|| \cdot ||_2$ denotes  the standard spectral norm of a matrix or a vector. 

\bigskip

The justification for the statistics used from above is given by the following theorem 
(for details see \cite{Nick}, pp.  132, 370-371):
\begin{theorem}\label{thm6}
Let  $A \in \mathbb R^{n \times n}$  and $0 \leq \epsilon<1$. Then
\begin{enumerate}
\item  $\|I_n- A^T A \|_2 \leq \epsilon$ $\Leftrightarrow$  There exists an orthogonal matrix $Q$ and  $E$ such that 
$A=Q+ E$,  where  $\|E\|_2 \leq \epsilon.$
That is, the matrix $A$ is very close to the true orthogonal matrix.
\item $\|A e-e\|_2  \leq \epsilon$ $\Leftrightarrow$  There exists $E_1$ such that $(A+E_1)e=e$, where  $\| E_1\|_2 \leq \frac{1}{\sqrt{n}}\,\epsilon$.
That is, all of $A+E_1$ row sums equal one.
\item $\|A^T e-e\|_2  \leq \epsilon$  $\Leftrightarrow$ There exists $E_2$ such that $(A+E_2)^Te=e$, where  $\| E_2\|_2 \leq \frac{1}{\sqrt{n}}\,\epsilon$.
That is, all of $A+E_2$ column sums equal  one.
 \end{enumerate}
\end{theorem}

\pagebreak

Several examples to test  our algorithms are considered.
\begin{example}\label{example1a}
We present a comparison of Algorithm $1$ and Algorithm $1a$ for $z$ very close to $0$. 
Notice that the matrices $A$ generated by these two methods for the same value of $z$ may be completely different.
We see that the catastrophic cancellation occurs
in  Algorithm 1a for $z \approx 0$, see Table $2$. In contrast, Algorithm 1 gives perfectly accurate results, see Table $1$.

\medskip

\begin{table}\label{tabelka1a}
\caption{The results for Example $1$  and  the  matrices $A(3 \times 3)$ computed by Algorithm $1$. }

\medskip
\begin{center}
\begin{tabular}{|c|c|c|c|c|c|}
\hline
$z$ &  $10^{-3}$   & $10^{-6}$  & $10^{-9}$  & $10^{-12}$ & $10^{-14}$\\
\hline ${err}_{orth}$ &   $9.12E-20$  & $2.22E-16$  & $5.65E-26$   & $4.84E-29$ & $6.03E-31$\\  
\hline ${err}_{rows}$ &   $1.11E-16$  & $2.71E-16$  & $0$   & $0$ & $0$\\  
\hline ${err}_{columns}$ &   $1.11E-16$  & $2.71E-16$  & $0$   & $0$ & $0$\\  
\hline
\end{tabular}
\end{center}
\end{table}

\medskip

\begin{table}
\label{tabelka1b}
\caption{The results for Example $1$  and  the  matrices $A(3 \times 3)$ computed by Algorithm $1a$. }

\medskip
\begin{center}
\begin{tabular}{|c|c|c|c|c|c|}
\hline
$z$ &  $10^{-3}$   & $10^{-6}$  & $10^{-9}$  & $10^{-12}$ & $10^{-14}$\\
\hline ${err}_{orth}$ &   $1.23E-13$  & $1.12E-11$  & $1.65E-07$   & $1.55E-04$ & $0.0016$\\  
\hline ${err}_{rows}$ &   $1.07E-13$  & $9.75E-12$  & $1.43E-07$   & $1.34E-04$ & $0.0014$\\  
\hline ${err}_{columns}$ &   $1.03E-13$  & $9.75E-12$  & $1.43E-07$   & $1.34E-04$ & $0.0014$\\  
\hline
\end{tabular}
\end{center}
\end{table}
\end{example}

\medskip

\begin{example}\label{example1}
We test  Algorithm $2$  on random matrices $X(n \times n)$ generated by the \textsl{MATLAB} code:

\begin{verbatim}
randn('state',0);
X=randn(n);
Q=Algorithm2(X);
err_orth=norm(eye(n)-Q'*Q);
\end{verbatim}

Random matrices of entries are  from the normal distribution ${\cal N}(0,1)$. They
are generated by the \textsl{MATLAB} function "randn". 
Before each call, the random number generator is reset to its
initial state.

Visibly Algorithm $2$ gives very satisfactory results (see Table \ref{tabelka1}). 
Theorem \ref{thm6} guarantees that every computed matrix $Q$ is very close to the exactly orthogonal matrix. 

\medskip

\begin{table}\label{tabelka1}
\caption{The orthogonality error  for Example $2$ and  the  matrix $Q(n \times n)$ computed by Algorithm $2$. }

\medskip

\begin{center}
\begin{tabular}{|c|c|c|c|c|c|}
\hline
$n$ &  $10$   & $50$  & $100$  & $500$ & $1000$\\
\hline ${err}_{orth}$ &   $1.34E-15$  & $2.21E-15$  & $2.30E-15$   & $3.40E-15$ & $6.34E-15$\\  
\hline
\end{tabular}
\end{center}
\end{table}
\end{example}

\medskip

\begin{example}\label{example2}
In the next step we test Algorithm $3$  on matrices $Q(n \times n)$ generated by Algorithm $2$ 
as specified in Example \ref{example1} and 
on orthogonal matrices $W((n-1) \times (n-1))$ generated by Householder QR decomposition of random matrices.

The following \textsl{MATLAB} code is used:

\begin{verbatim}
randn('state',0);
X=randn(n); Q=Algorithm2(X);
Y=randn(n-1); [W,R]=qr(Y);
A=Algorithm3(Q,W);
e=ones(n,1);
\end{verbatim}

Again, as illustrated in Table \ref{tabelka2}, Algorithm $3$ yields  very good results. 

\medskip

\begin{table}\label{tabelka2}
\caption{The results for Example $3$ and  the  matrix $A(n \times n)$ computed by Algorithm $3$. }

\medskip
\begin{center}
\begin{tabular}{|c|c|c|c|c|c|}
\hline
$n$ &  $10$   & $50$  & $100$  & $500$ & $1000$\\
\hline ${err}_{orth}$ &   $1.34E-15$  & $2.21E-15$  & $2.30E-15$   & $3.40E-15$ & $6.34E-15$\\  
\hline ${err}_{rows}$ &   $1.12E-15$  & $4.50E-15$  & $7.26E-15$   & $2.43E-14$ & $4.34E-14$\\  
\hline ${err}_{columns}$ &   $1.09E-15$  & $4.69E-15$  & $7.56E-15$   & $2.39E-14$ & $4.03E-14$\\  
\hline
\end{tabular}
\end{center}
\end{table}
\end{example}

\medskip

\begin{example}\label{example3}
We test now Algorithm $4$ with  the following  \textsl{MATLAB} code:
\begin{verbatim}
randn('state',0);
i=sqrt(-1); r=2;p=3;z=[0.6+0.8*i,-0.8+0.6*i];
n=r+p+4;
X=randn(n); Q=Algorithm2(X);
A=Algorithm4(r,p,z,Q);
eigA=eig(A) % The vector eigA contains the computed  eigenvalues of A
\end{verbatim}
The exact eigenvalues of $A$  are: $1,1,-1,-1,-1,0.6 \pm0.8 i, -0.8 \pm 0.6$.\\

The corresponding eigenvalues of computed matrix $A$ generated by Algorithm $4$ are:
\begin{verbatim}
eigA =
      6.000000000000001e-01 + 8.000000000000002e-01i
      6.000000000000001e-01 - 8.000000000000002e-01i
     -7.999999999999996e-01 + 5.999999999999996e-01i
     -7.999999999999996e-01 - 5.999999999999996e-01i
      1.000000000000000e+00 + 0.000000000000000e+00i
      1.000000000000000e+00 + 0.000000000000000e+00i
     -1.000000000000000e+00 + 0.000000000000000e+00i
     -1.000000000000000e+00 + 0.000000000000000e+00i
     -9.999999999999998e-01 + 0.000000000000000e+00i
\end{verbatim}

Thus, upon comparing the latter, the statistics
\[
err_{orth} = 1.08E-15, \quad err_{rows} = 8.88E-16, \quad err_{columns} =9.15E-16,
\]
renders all results almost perfect in floating-point arithmetic. 
\end{example}

\medskip

\begin{example}\label{example4}
Finally, the performance of Algorithm $6$  is tested. In doing so, the following \textsl{MATLAB} code is used:
\begin{verbatim}
n=2; m=n^2;d=[1,-1,1,1];
B=Algorithm5(n,d); 
randn('state',0); X=randn(n); P=Algorithm2(X); 
A=Algorithm6(B,P)
e=ones(m,1);
err_orth=norm(eye(m)-A'*A)
err_rows=norm(A*e-e)
err_columns=norm(A'*e-e)
\end{verbatim}

The outcoming statistics read as:
\[
err_{orth} = 8.55E-16, \quad err_{rows} = 9.28E-16, \quad err_{columns} =9.15E-16.
\]

Clearly all results produces high accuracy in floating-point arithmetic's.
 Recall that in the first step of Algorithm $6$, the Algorithm $5$ is applied.
\end{example}

\bigskip

\section{Appendix - \textsl{MATLAB} Codes}

For the sake of completeness,  we enclose \textsl{MATLAB} codes to all discussed Algorithms in question.

\begin{verbatim}
function  [A]=Algorithm1(z)  
% [A]=Algorithm1(z) 
% A(3x3) is orthogonal and symmetric generalized stochastic matrix.
% Parameter z should be in the interval [-1/3,1].
n=3; A=zeros(n);
if z>1 || z<-1/3
   disp('z should be in the interval [-1/3,1]');
   return;
end
t=1-z; 
if t==0
   x=0; y=0;
   A=[0 0 1;0 1 0;1 0 0];
   return;
end
delta=t*(1+3*z);
x=(t+sqrt(delta))/2;
y=-z*t/x(1);
A=[x y z;y z x;z x y];
end
\end{verbatim}

\bigskip
\begin{verbatim}
function  [A]=Algorithm1a(z)  
% [A]=Algorithm1(z) (unstable for z  close to 0)
% A(3x3) is orthogonal and symmetric generalized stochastic matrix.
% Parameter z should be in the interval [-1/3,1].
n=3; A=zeros(n);
if z>1 || z<-1/3
   disp('z should be in the interval [-1/3,1]');
   return;
end
t=1-z; 
if t==0
   A=[0 0 1;0 1 0;1 0 0];
   return;
end
delta=t*(1+3*z);
if z==0
   A=[0 1 0;1 0 0;0 0 1];
return;
end
x=(t-sqrt(delta))/2;
y=-z*t/x;
A=[x y z;y z x;z x y];
end
\end{verbatim}
 
\bigskip

\begin{verbatim}
function  [Q]=Algorithm2(X)  
% [Q]=Algorithm2(X).  
% Q(nxn) is orthogonal and g.d.s.
% The first column of Q is e/sqrt(n), where e=(1,1,...,1).
% Householder Q-R decomposition is used.
[m,n]=size(X);
Q=zeros(n);
if m~=n
   disp('X should be a square matrix.');
   return;
end
e=ones(n,1);
norm_e=sqrt(n);
X(:,1)=e/norm_e;
[Q,~]=qr(X);  Q=-Q;
end
\end{verbatim}

\bigskip

\begin{verbatim}
function  [A]=Algorithm3(Q,W)  
% [A]=Algorithm3(Q,W) 
% A(nxn)is orthogonal generalized doubly stochastic matrix.
% Q(nxn) is an orthogonal matrix with the first column e/sqrt(n).
% W(n-1)x(n-1) is an orthogonal matrix. 
[m,n]=size(Q);
A=zeros(n);
if m~=n
   disp('X should be a square matrix.');
   return;
end
[k,l]=size(W);
if k~=l
   disp('Y should be a square matrix.');
   return;
end
if k~=(n-1)
   disp('Size of Y should be equal to n-1');
   return;
end
z=zeros(n-1,1);
B=[1,z'; z,W];
A=Q*B*Q';
end
\end{verbatim}

\bigskip

\begin{verbatim}
function  A=Algorithm4(r,p,z,Q)  
%[A]=Algorithm4(r,p,z,Q)  
% A(nxn) is orthogonal  generalized doubly stochastic.
% n=r+p+2m, where m is the length of a vector z,
% r is the number of 1's, and p is the number of -1's of A.
% Here |z(k)|=1 for k=1,..., m.
% Assume that r>=1,  p>=1, and m>=1. 
% Q(nxn) is an orthogonal matrix with the first column e/sqrt(n).
z=z(:);
m=length(z);
n=r+p+2*m;
A=eye(n);
c=real(z);s=imag(z); 
R=zeros(2*m,2*m);
for k=1:m
    Rk=[c(k) s(k);-s(k) c(k)];
    R(2*k-1:2*k,2*k-1:2*k)=Rk;
end
B=[eye(r) zeros(r,p) zeros(r,2*m)
   zeros(p,r) -eye(p) zeros(p,2*m)
   zeros(2*m,r)  zeros(2*m,p) R];
A=Q*B*Q';
end
\end{verbatim}

\bigskip

\begin{verbatim}
function [X] = Algorithm5(n,d)
% [X] = Algorithm5(n,d)
% X(mxm), m=n^2, X is a solution of the YBE
% d=(d(1),..., d(m)), where m=n^2.
m=max(size(d));
if ~(m==n*n)
    disp('Wrong dimensions');
return;
end
for j=1:n
    for i=1:n
        S(i,j)=(i-1)*n+j;
    end
end
p=[];
for i=1:n
    p=[p; S(:,i)];
end
p=p'; X=diag(d); X=X(:,p);
end
\end{verbatim}

\bigskip

\begin{verbatim}
function  [A]=Algorithm6(B,P)  
% [A]=Algorithm6(B,P) 
% B(mxm), m=n^2, satisfies the Yang-Baxter equation.
% Be1=e1 and B'e1=e1, where e1=(1,0,...,0)'.
% P(nxn) is orthogonal with the first column e/sqrt(n), where e=(1,1,...,1)'.
% A(mxm) is orthogonal generalized doubly stochastic matrix satisfying the Yang-Baxter equation.
[m,m]=size(B);
[n,n]=size(P);
A=eye(m);
if m~=n^2
   disp('Wrong dimensions!');
   return;
end
Q=kron(P,P);
A=Q*B*Q';
end
\end{verbatim}
\end{document}